\pdfoutput=1
\documentclass[twoside]{article}
\usepackage{aistats2012mod}

\usepackage{amsmath}
\usepackage{amssymb}
\usepackage{amsthm}
\usepackage{amsfonts,graphicx,color}
\usepackage{subfigure}
\usepackage{fullpage}
\usepackage{framed}
\usepackage{times}
\usepackage[sort]{natbib}
\newtheorem{thm}{Theorem}

\newtheorem{prop}{Proposition}

\def \E{\mathbb{E}}

\def \cE{{\cal E}}
\def \bomega{\mathbf{\omega}}
\newcommand{\Leb}{\mathop{\mathrm{Leb}}}
\newcommand{\Vol}{\mathop{\mathrm{Vol}}}

\newcommand{\diam}{\mathop{\mathrm{diam}}}

\newcommand{\vol}{\mathop{\mathrm{Vol}}}
\newcommand{\leb}{\mathop{\mathrm{Leb}}}

\let\hat\widehat

\begin{document}

%

%

\twocolumn[

\aistatstitle{Adaptive Semisupervised Inference}

\aistatsauthor{ Martin Azizyan \And Aarti Singh \And Larry Wasserman }

\aistatsaddress{ 
Machine Learning Department \\
Carnegie Mellon University \\
Pittsburgh, PA, 15213-3890 USA
 \And 
Machine Learning Department\\
Carnegie Mellon University\\
Pittsburgh, PA, 15213-3890 USA
 \And
Department of Statistics and\\
Machine Learning Department\\
Carnegie Mellon University\\
Pittsburgh, PA, 15213-3890 USA } ]

\begin{abstract}
Semisupervised methods inevitably invoke some assumption that links the marginal 
distribution of the features to the regression function of the label.
Most commonly, the cluster or manifold assumptions are used which imply that
the regression function is smooth over high-density clusters or manifolds supporting the
data.
A generalization of these assumptions is that the regression function is smooth with 
respect to some density sensitive distance. This motivates the use of a density 
based metric \citep{bousquet04,orlitsky,diff_maps} for semisupervised learning. 
We analyze this setting and make the following contributions  - 
(a) we propose a semi-supervised learner that uses a density-sensitive kernel
and show that it provides better performance than any supervised learner if
the density support set has a small condition number 
and 
(b) we show that it is possible to adapt to the 
degree of semi-supervisedness using data-dependent choice of a parameter
that controls sensitivity of the distance metric to the density. This ensures that the
semisupervised learner never performs worse than a supervised learner even if the
assumptions fail to hold.
\end{abstract}

\section{Introduction}

Semisupervised methods inevitably invoke some assumption that links
the marginal distribution $p(x)$ of the features $X$ to the
regression function $f(x) = \E[Y|X = x]$ of the label $Y$. The most
common assumption is the {\em cluster assumption} in which it is
assumed that $f$ is very smooth wherever $p$ exhibits clusters
\citep{seeger,rigollet07,LW:nips07,ASingh:unlabeled}.  In the special
case where the clusters are manifolds, this is called the {\em
manifold assumption} \citep{LW:nips07,belkin_niyogi,partha}.

A generalization of the cluster and manifold assumptions is that the
regression function is smooth with respect to some density-sensitive
distance. Several recent papers propose using a density based metric or
diffusion distance for semisupervised learning
\citep{orlitsky,diff_maps,bousquet04}. In this paper, we analyze semisupervised
inference under this generalized assumption.

Singh, Nowak and Zhu \citeyearpar{ASingh:unlabeled}, Lafferty and
Wasserman \citeyearpar{LW:nips07} and
Nadler et al \citeyearpar{nadler09}
have showed that the degree to which
unlabeled data improves performance is very sensitive to the cluster
and manifold assumptions. 
In this paper, we introduce {\em adaptive semisupervised
inference}. We define a parameter $\alpha$ that controls the
sensitivity of the distance metric to the density, and hence the
strength of the semisupervised assumption. When $\alpha = 0$ there is
no semisupervised assumption, that is, there is no link between $f$
and $p$. When $\alpha = \infty$ there is a very strong semisupervised
assumption. We use the data to estimate $\alpha$ and hence we adapt to
the appropriate assumption linking $f$ and $p$.


This paper makes the following contributions - (a) we propose a
semi-supervised learner that uses a density-sensitive kernel and show
that it provides better performance than any supervised learner 
if the density support set has a small condition number 
and (b) we show that it is possible to adapt to the degree of
semi-supervisedness using data-dependent choice of a parameter that
controls sensitivity of the distance metric to the density. This
ensures that the semisupervised learner never performs worse than a
supervised learner even if the assumptions fail to hold.
Preliminary simulations, to be reported in future work, confirmed that our proposed estimator adapts
well to alpha and has good risk when the semisupervised smoothness
holds and when it fails.  

{\em Related Work.}
There are a number of papers that
discuss conditions under which 
semisupervised methods can succeed
or that discuss metrics that are useful for
semisupervised methods.
These include
\citet{bousquet04}, \citet{SSL_TR},
\citet{nadler09}, \citet{orlitsky}
and references therein.
However, to the best of our knowledge,
there are no papers that explicitly study
adaptive methods that allow the data to choose the
strength of the semisupervised assumption.

{\em Outline.}
This paper is organized as follows. In Section~\ref{sec:setup} we
define a set of joint distributions ${\cal P}_{XY}(\alpha)$ indexed
by $\alpha$. In Section~\ref{sec:est}, we define a density sensitive
estimator $\widehat f_\alpha$ of $f$, assuming that $(f, p) \in {\cal
P}_{XY}(\alpha)$. We find finite sample bounds on the error of
$\widehat f_\alpha$ and we investigate the dependence of this error on
$\alpha$.  In Section~\ref{sec:adap}, we show that
cross-validation can be used to adapt to $\alpha$.
We conclude in
Section~\ref{sec:disc}.

\section{Definitions}
\label{sec:setup}

We consider the collection of joint distributions ${\cal P}_{XY}
(\alpha) = {\cal P}_X \times {\cal P}_{Y|X}$ indexed by a
density-sensitivity parameter $\alpha$ as follows. $X, Y$ are random
variables, $X$ is supported on a compact domain ${\cal X} \subset \mathbb{R}^d$, and $Y$ is real-valued.
The marginal density $p(x) \in [\lambda_0,\Lambda_0]$ is bounded over its support 
$\{x:p(x)>0\}$, where $0<\lambda_0,\Lambda_0<\infty$.
Also, let the conditional density be $p(y|x)$ with variance bounded by $\sigma^2$,
and conditional label mean or regression function be 
$f(x)= \E[Y|X=x]$, with $|f(x)|\leq M$. 
We say that $(p, f) \in {\cal P}_{XY}(\alpha)$ if these
functions satisfy the properties described below.

Before stating the properties of $f$ and $p$, we define a distance
metric with density sensitivity $\alpha$.

{\bf Density-sensitive distance:} We consider the following distance
with density sensitivity $\alpha \in [0,\infty)$ between two points
$x_1, x_2 \in {\cal X}$ that is a modification of the definition in
\citet{orlitsky}:
\begin{equation}
D_\alpha(x_1,x_2) = 
\inf\limits_{\gamma\in\Gamma(x_1,x_2)}
\int\limits_0^{L(\gamma)} \frac1{p(\gamma(t))^\alpha} dt,
\end{equation}
where $\Gamma(x_1,x_2)$ is the set of all continuous finite curves
from $x_1$ to $x_2$ with unit speed everywhere and $L(\gamma)$ is the
length of curve $\gamma$ (i.e. $\gamma(L(\gamma))=x_2$).
Notice that 
large $\alpha$ makes points connected by high density paths
closer, and $\alpha=0$ corresponds to Euclidean distance.

Our first assumption is that
the regression function $f$ is smooth with respect to the density sensitive distance:

{\bf A1) Semisupervised smoothness:} The regression
function $f(x) = \E[Y|X=x]$ is $\beta$-smooth with respect to the
density-sensitive distance $D_\alpha$, i.e. there exists constants
$C_1,\beta>0$ such that for all $x_1,x_2 \in {\cal X}$
$$
|f(x_1) - f(x_2)| \leq C_1 \ \Bigl[D_\alpha(x_1,x_2)\Bigr]^\beta.
$$
In particular if $\alpha=0$ and $\beta=1$, this corresponds to Lipschitz smoothness.

Our second assumption is that the density function $p$ is smooth 
with respect to Euclidean distance over the support set.
Recall that the {\em support} of $p$ is
$S = \{x:\ p(x) > 0\}$.

{\bf A2) Density smoothness:} The density function $p(x)$ is H\"{o}lder
$\eta$-smooth with respect to Euclidean distance if it has $\lfloor \eta \rfloor$
derivatives and there exists a
constant $C_2>0$ such that for all $x_1,x_2 \in S$
$$
|p(x_1) - T^{\lfloor \eta \rfloor}_{x_2}(x_1)| \leq C_2 \ \|x_1-x_2\|^\eta,
$$
where $\lfloor \eta \rfloor$ is the largest integer such that $\lfloor \eta \rfloor < \eta$, 
and $T^{\lfloor \eta \rfloor}_{x_2}$ is the Taylor polynomial of degree $\lfloor \eta \rfloor$ around
the point $x_2$.

The {\em condition number} of a set $S$ with boundary $\partial S$
is the largest real number $\tau>0$ such that,
if $d(x,\partial S) \leq \tau$ then $x$ has a unique projection onto the boundary of $S$.
Here, $d(x,\partial S) = \inf_{z\in \partial S}||x-z||$.
When $\tau$ is large, $S$ cannot be too thin, the boundaries of $S$ cannot be too curved
and $S$ cannot get too close to being self-intersecting.
If $S$ consists of more than one connected component,
then $\tau$ large also means that the connected components cannot be too close to each other.
Let $\tau_0$
denote the smallest condition number of the support sets $S$ of all $p\in\mathcal{P}_X$.
We shall see that semisupervised inference outperforms supervised inference
when $\tau_0$ is small. Additionally, we assume that $S$ has at most $K<\infty$
connected components.

In the supervised setting, we assume access to $n$ labeled data 
${\cal L} = \{X_i, Y_i\}^n_{i=1}$ drawn i.i.d. from 
${\cal P}_{XY}(\alpha)$, and in the semi-supervised setting, we assume access to $m$
additional unlabeled data ${\cal U} = \{X_i\}^m_{i=1}$ drawn
i.i.d. from ${\cal P}_{X}$.

As usual, we write
$a_n = O(b_n)$ if
$|a_n/b_n|$ is bounded for all large $n$.
Similarly,
$a_n = \Omega(b_n)$ if
$|a_n/b_n|$ is bounded away from 0 for all large $n$.
We write
$a_n \asymp b_n$ if
$a_n = O(a_n)$ and
$a_n = \Omega(b_n)$.

\section{Density-Sensitive Inference}
\label{sec:est}

Let $K(x)$ be a symmetric non-negative function and let
$K_h(x) = K(\|x\|/h)$.
Let
\begin{equation}
\hat p_m(x) = \frac{1}{m}\sum_{i=1}^m \frac{1}{h_m^d} K_{h_m}(x-X_i)
\end{equation}
be the kernel density estimator of $p$ with bandwidth $h_m$,
based on the unlabeled data.
Define the support set estimate $\hat S = \{x:\hat p_m(x) > 0\}$ and the empirical
boundary region
$$
\widehat{\mathcal{R}}_{\hat\partial S} = \left\{x:\inf\limits_{z\in\partial \widehat{S}} \|x-z\|_2 < 2\delta_m\right\}.
$$
where $\delta_m=2c_2\sqrt{d}\left((\log^2m)/m\right)^{\frac{1}{d}}$ for some constant $c_2>0$.
Now define a plug-in estimate of the $D_\alpha$ distance as follows:
$$
\hat D_{\alpha,m}(x_1,x_2) = 
\inf\limits_{\gamma\in\hat\Gamma(x_1,x_2)}\int\limits_0^{L(\gamma)} 
\frac1{\hat p_m(\gamma(t))^\alpha} dt,
$$ 
where $\hat \Gamma(x_1,x_2) = \{\gamma \in \Gamma(x_1,x_2): \forall t \in
[0,L(\gamma)] \ \gamma(t) \in \hat S \setminus \hat{\mathcal{R}}_{\hat \partial S}\}$, and $\hat D_{\alpha,m}(x_1,x_2) = \infty$
if $\hat \Gamma(x_1,x_2) = \emptyset$.

We consider the following 
semisupervised
learner which uses a kernel that is
sensitive to the density.
In the following definitions we take, for simplicity,
$K(x) = I(||x|| \leq 1)$.


{\bf Semisupervised kernel estimator:}
\begin{equation}
\widehat f_{h,\alpha}(x) = \frac{\sum^n_{i=1}Y_i K_h\left(\widehat D_{\alpha,m}(x,X_i)\right)}
{\sum^n_{i=1}K_h\left(\widehat D_{\alpha,m}(x,X_i)\right)}.
\end{equation}

\subsection{Performance upper bound for semisupervised estimator}
\label{sec:SSL_UB}
The following theorem characterizes the performance of the density sensitive
semisupervised kernel estimator.
\begin{thm}
\label{thm:SSL_UB}
Assume $\lambda_0>1+c_0$ for some constant $c_0>0$
\footnote{This assumption is more restrictive than necessary, and a more general statement can be by introducing a rescaling factor in the definition of the density-sensitive distance.}
and let $\epsilon_m=c_1 (\log m)^{-1/2}$ for constant 
$c_1 >0$ 
and $\delta_m=2c_2\sqrt{d}\left((\log^2m)/m\right)^{\frac{1}{d}}$ 
for some constant $c_2>0$.
If  $\tau_0\in(3\delta_m, \infty)$ and
$h > (2c_4/(\tau_0^{d-1}(\lambda_0-\epsilon_m)^\alpha))$
where $c_4>0$ is a constant,
then for large enough $m$ 
\begin{align*}
& \sup\limits_{(p,f)\in\mathcal{P}_{XY}(\alpha)} \mathbb{E}_{n,m}\left\{\int (\widehat{f}_{h,\alpha}(x)-f(x))^2 dP(x)\right\} \leq\\
& \hspace{3cm}
(M^2+\sigma^2)\left(\frac{1}{m}+3 c_3 2^d\Lambda_0 \frac{\delta_m}{\tau_0}\right) \\
&\hspace{3cm}+ \left[h \left(\frac{\lambda_0+\epsilon_m}{\lambda_0}\right)^\alpha \right]^{2\beta}\\
&\hspace{3cm}+ \frac{K(M^2/e+2\sigma^2)}{n}.
\end{align*}
\end{thm}

The proof of Theorem~\ref{thm:SSL_UB} is given in section~\ref{sec:SSL_UB_proof}.
The first term is negligible when the amount of unlabeled data $m$ is large.
The second term is the bias and third term is variance. 
If the bandwidth
$$h \asymp \frac1{\delta_m^{d-1}\lambda_0^\alpha}$$ 
and $\alpha \asymp \log m$ is large enough, then the density-sensitive 
semisupervised kernel estimator 
is able to achieve an integrated MSE rate of $O(n^{-1})$ 
for all joint distributions in ${\cal P}_{XY} (\alpha)$ supported on sets with condition number
$\tau_0 >3\delta_m$.   

\subsection{Performance lower bound for any supervised estimator}
We now establish a lower bound on the performance of any supervised estimator.
\begin{thm}
Assume $d\geq2$ and $\alpha>0$.
There exists a constant $c_5>0$ depending only on $d$ so that if $\tau_0\leq c_5 n^{-\frac{1}{d-1}}$, then
\begin{align*}
\inf_{\widehat{f}}\sup_{(p,f)\in\mathcal{P}_{XY}(\alpha)} \mathbb{E}_n\int(\widehat{f}(x)-f(x))^2 dP(x) = \Omega(1)
\end{align*}
where the inf is over all supervised estimators.
\label{thm:SL_LB}
\end{thm}
Coupled with Theorem~\ref{thm:SSL_UB}, the results state that if the condition number of 
the support set is small $3\delta_m < \tau_0 \leq c_5 n^{-\frac1{d-1}}$ and $\alpha$ is large enough, then 
the density-sensitive semi-supervised estimator outperforms any supervised learning 
algorithm in terms of integrated MSE rate. 

A complete proof of Theorem~\ref{thm:SL_LB} is given in the appendix.
Here we provide some intuition regarding the proof strategy.
We construct a set of joint distributions over $X$ and $Y$ that depends on $n$, and apply Assouad's Lemma.
Intuitively, we need to take advantage of the decreasing condition number $\tau_0$.
This is because if $\tau_0$ were to be kept fixed, as $n$ increases the semi-supervised assumption would reduce to familiar Euclidean smoothness.

\begin{figure}
  \centering
    \includegraphics[width=0.49\textwidth,clip=true,trim=6.6cm 4cm 6.5cm 3.3cm]{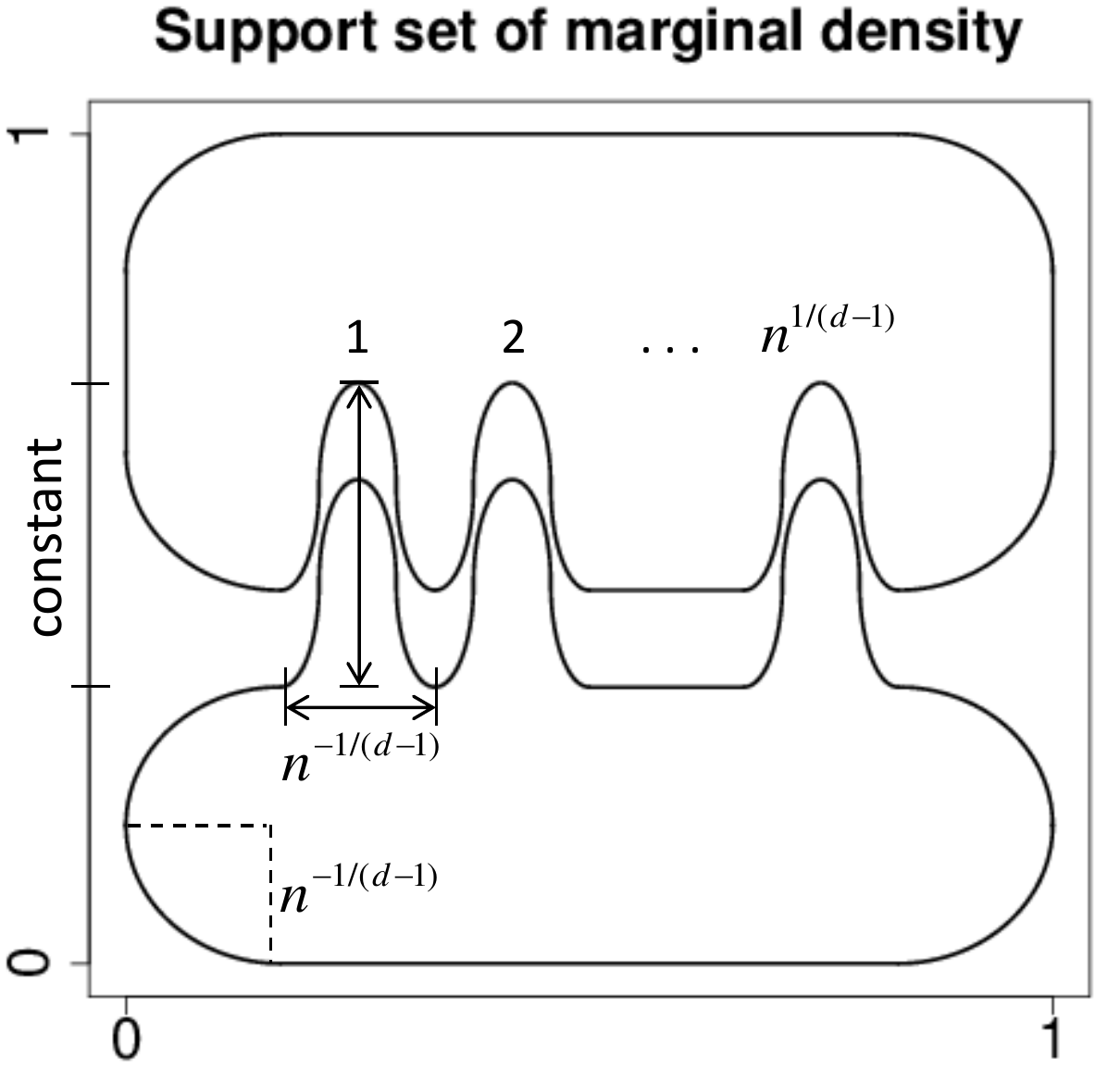}
\vspace{-25pt}
  \caption{A two-dimensional cross-section of the support of a marginal density $p$ used in the proof of Theorem~\ref{thm:SL_LB}.}
  \label{fig:lb}
\end{figure}

So, we construct the distributions as follows.
We split the unit cube in $\mathbb{R}^d$ into two rectangle sets with a small gap in between, and let the marginal density $p$ be uniform over these sets.
Then we add a series of ``bumps'' between the two rectangles, as shown schematically in Figure \ref{fig:lb}.
Over one of the sets we set $f\equiv M$, and over the other we set $f\equiv -M$.
The number of bumps increases with $n$, implying that the condition number must decrease.
The sets are designed specifically so that the condition number can be lower bounded easily as a function of $n$.
In essence, as $n$ increases these boundaries become space-filling, so that there is a region where the regression function could be $M$ or $-M$, and it is not possible to tell which with only labeled data.

\section{Adaptive Semisupervised Inference}
\label{sec:adap}

In section~\ref{sec:SSL_UB}, we established a bound on the integrated mean
square error of the density-sensitive semisupervised kernel estimator.
The bound is achieved by using an estimate $\hat D_\alpha$ of the 
density-sensitive distance. 
However, this requires knowing the density-sensitive parameter $\alpha$, along with
other parameters.

It is critical to choose $\alpha$ (and $h$) appropriately, otherwise we might
incur a large error if the semisupervised assumption does not hold or
holds with a different density sensitivity value $\alpha$. The
following result shows that we can adapt to the correct degree of
semisupervisedness $\alpha$ if cross-validation is used to select the
appropriate $\alpha$ and $h$. This implies that the estimator
gracefully degrades to a supervised learner if the semisupervised
assumption (sensitivity of regression function to marginal density)
does not hold ($\alpha = 0$).

For any $f$, define the risk
$R(f) = \E[(f(X)-Y)^2]$ and the excess risk
$\cE(f) = R(f) - R(f^*) = \E[(f(X)-f^*(X))^2]$
where $f^*$ is the true regression function.
Let
${\cal H}$ be a finite set of bandwidths
and let
${\cal A}$ be a finite set of values for $\alpha$.
Divide the data into training data $T$
and validation data $V$.
For notational simplicity, let both sets have size $n$.
Let 
${\cal F} = \{\hat f^T_{\alpha,h}\}_{\alpha \in {\cal A}, h\in {\cal H}}$ 
denote the semisupervised kernel estimators trained on
data $T$ using $\alpha \in {\cal A}$ and $h \in {\cal H}$.  
For each $\hat f_{\alpha,h}^T\in {\cal F}$ let
$\hat R^V (\hat f^T_{\alpha,h}) = n^{-1}\sum^n_{i=1}(\hat f^T_{\alpha,h}(X_i)-Y_i)^2$ 
where the sum is over $V$.
Let $Y_i = f(X_i) + \epsilon_i$ with $\epsilon_i \stackrel{i.i.d}{\sim} {\cal N}(0,\sigma^2)$.
Also, we assume that 
$|f(x)|, |\hat f^T_{\alpha,h}(x)| \leq M$, where $M>0$ is a constant.\footnote{
Note that the estimator can always be truncated if necessary.}

\begin{thm}
\label{thm:crossval}
Let ${\cal F} = \{\hat f^T_{\alpha,h}\}_{\alpha \in {\cal A}, h\in
{\cal H}}$ denote the semisupervised kernel estimators trained on
data $T$ using $\alpha \in {\cal A}$ and $h \in {\cal H}$.  Use
validation data $V$ to pick
$$
(\hat \alpha,\hat h)  = 
\arg\min_{(\alpha \in {\cal A},h\in {\cal H})} \hat R^V(\hat f^T_{\alpha,h})
$$
and define the corresponding estimator $\hat f_{\hat \alpha,\hat h}$. Then,
for every $0 < \delta < 1$,
\begin{align*}
\E[\cE(\hat f_{\hat \alpha,\hat h})] \leq \frac1{1-a}
&\left[\min_{\alpha\in {\cal A}, h \in {\cal H}} \E[\cE(\hat
  f_{\alpha,h})]\right.\\ 
&\hspace{0.7cm}+ \left.\frac{\log(|{\cal A}||{\cal H}|/\delta)}{nt}
\right] + 4\delta M^2
\end{align*}
where $0<a<1$ and $0<t <  15/(38(M^2+\sigma^2))$ are constants. $\E$ denotes expectation over 
everything that is random.
\end{thm}

See appendix for proof.
In practice, both
${\cal H}$ and ${\cal A}$ may be taken to be 
of size $n^a$ for some $a>0$.
Then we can approximate the optimal $h$ and $\alpha$
with sufficient accuracy to achieve the optimal rate.
Setting $\delta = 1/(4 M^2n)$, we then see that
the penalty for 
adaptation is
$\frac{\log(|{\cal A}||{\cal H}|/\delta)}{nt} + \delta M = O(\log n /n)$
and hence introduces only a logarithmic term.



\section{Discussion}
\label{sec:disc}

Semisupervised methods are very powerful
but, like all methods, they only work under
certain conditions.

We have shown that, when the support of the distribution is somewhat 
irregular
(i.e., the boundary of the support of the density has a small condition 
number),
then semi-supervised methods can attain better performance.
Specifically, we demonstrated that a semi-supervised kernel estimator that
uses a density-sensitive distance can outperform any supervised 
estimator in
such cases.

We introduced a family of estimators
indexed by a parameter $\alpha$.
This parameter controls the strength of the semi-supervised
assumption.
We showed that the behavior of the semi-supervised method
depends critically on $\alpha$.

Finally, we showed that cross-validation
can be used to automatically adapt to $\alpha$ so that $\alpha$ does not 
need to be known.
Hence, our method takes advantage of the unlabeled data
when the semi-supervised assumption holds, but does not add extra
bias when the assumption fails.
Preliminary simulations confirm that our proposed estimator adapts
well to alpha and has good risk when the semi-supervised smoothness
holds and when it fails.  We will report these results in future work.

The analysis in this paper can be extended in several ways.
First, it is possible to use other density sensitive metrics
such as the diffusion distance \citep{wasserman08spectral}.
Second, it is possible to relax the assumption that
the density $p$ is strictly bounded away from 0 on its support.
Finally, other estimators besides kernel estimators can be used.
We will report on these extensions elsewhere.

\section{Proof of Theorem~\ref{thm:SSL_UB}}
\label{sec:SSL_UB_proof}
Here we prove Theorem~\ref{thm:SSL_UB} stated in section~\ref{sec:SSL_UB} (repeated below for convenience), using some results given in the appendix.

\begin{thm}
Assume $\lambda_0>1+c_0$ for some constant $c_0>0$
\footnote{This assumption is more restrictive than necessary, and a more general statement can be by introducing a rescaling factor in the definition of the density-sensitive distance.}
and let $\epsilon_m=c_1 (\log m)^{-1/2}$ for constant 
$c_1 >0$ 
and $\delta_m=2c_2\sqrt{d}\left((\log^2m)/m\right)^{\frac{1}{d}}$ 
for some constant $c_2>0$.
If  $\tau_0\in(3\delta_m, \infty)$ and
$h > (2c_4/(\tau_0^{d-1}(\lambda_0-\epsilon_m)^\alpha))$
where $c_4>0$ is a constant,
then for large enough $m$ 
\begin{align*}
& \sup\limits_{(p,f)\in\mathcal{P}_{XY}(\alpha)} \mathbb{E}_{n,m}\left\{\int (\widehat{f}_{h,\alpha}(x)-f(x))^2 dP(x)\right\} \leq\\
& \hspace{3cm}
(M^2+\sigma^2)\left(\frac{1}{m}+3 c_3 2^d\Lambda_0 \frac{\delta_m}{\tau_0}\right) \\
&\hspace{3cm}+ \left[h \left(\frac{\lambda_0+\epsilon_m}{\lambda_0}\right)^\alpha \right]^{2\beta}\\
&\hspace{3cm}+ \frac{K(M^2/e+2\sigma^2)}{n}.
\end{align*}
\end{thm}

\begin{proof}
Let $\mathcal{G}_m$ be the indicator of the event when the unlabeled sample is such that $\sup\limits_{x\in S\backslash\mathcal{R}_{\partial S}}|p(x)-\widehat{p}_m(x)|\leq\epsilon_m$ and $\partial\widehat{S}\subset\mathcal{R}_{\partial S}$.
From Theorem \ref{thm:densityest},
\begin{align*}
\mathbb{E}_{n,m}\left\{(1-\mathcal{G}_m)\int (\widehat{f}_{h,\alpha}(x)-f(x))^2 dP(x)\right\} \\\leq \frac{1}{m}(M^2+\sigma^2).
\end{align*}
We can write
\begin{align*}
&\mathbb{E}_{n,m}\left\{\mathcal{G}_m\int (\widehat{f}_{h,\alpha}(x)-f(x))^2 dP(x)\right\} \\ &=
\mathbb{E}_{n,m}\left\{\mathcal{G}_m\int_{S_m^*} (\widehat{f}_{h,\alpha}(x)-f(x))^2 dP(x)\right\} \\
&+ \mathbb{E}_{n,m}\left\{\mathcal{G}_m\int_{S\backslash S_m^*} (\widehat{f}_{h,\alpha}(x)-f(x))^2 dP(x)\right\} 
\end{align*}
where $S_m^*$ as defined in Proposition \ref{thm:dbdest2}.
For the boundary region we have
\begin{align*}
&\mathbb{E}_{n,m}\left\{\mathcal{G}_m\int\limits_{S\backslash S_m^*} (\widehat{f}_{h,\alpha}(x)-f(x))^2 dP(x)\right\} \\
&\leq (M^2+\sigma^2) P(S\backslash S_m^*)\\
&\leq \Lambda_0(M^2+\sigma^2) \Leb(S\backslash S_m^*)
\end{align*}
where $\Leb$ denotes the Lebesgue measure.
Since the radius of curvature of $\partial S$ is at least $\tau_0$, and $\tau_0>3\delta_m$, we have by Proposition \ref{thm:condnumarea},
\begin{align*}
\Leb(S\backslash S_m^*) &\leq \Vol(\partial S) \frac{\left(\tau_0+3\delta_m\right)^d-\tau_0^d}{\tau_0^{d-1}}\\
&\leq c_3\left[\left(1+\frac{3\delta_m}{\tau_0}\right)^d-1\right]\\
&\leq c_3 \sum\limits_{i=1}^{d} \binom{d}{i} \frac{3\delta_m}{\tau_0}\\
&\leq 3 c_3 2^d \frac{\delta_m}{\tau_0}
\end{align*}
where $\Vol$ denotes the $d-1$-dimensional volume on $\partial S$.
So
\begin{align*}
\mathbb{E}_{n,m}\left\{\mathcal{G}_m\int\limits_{S\backslash S_m^*} (\widehat{f}_{h,\alpha}(x)-f(x))^2 dP(x)\right\} \\
\leq 3 c_3 2^d\Lambda_0(M^2+\sigma^2) \frac{\delta_m}{\tau_0}.
\end{align*}

Following the derivation in Chapter 5 of \citet{gyorfi2002nonparametric}, we have
\begin{align*}
&\mathbb{E}_{n}\left\{\mathcal{G}_m\int\limits_{S_m^*} (\widehat{f}_{h,\alpha}(x)-f(x))^2 dP(x)\right\} \\&\leq
\mathcal{G}_m C_1^2 \sup\limits_{x\in S_m^*} \sup\limits_{x'\in S\cap S_{x,h}^{\widehat{D}_{\alpha,m}}} D_\alpha(x,x')^{2\beta}\\
&+\mathcal{G}_m\frac{M^2/e+2\sigma^2}{n} \mathcal{N}\left(S_m^*,\widehat{D}_{\alpha,m},\frac{h}{2}\right)
\end{align*}
where $S_{x,h}^{\widehat{D}_{\alpha,m}}=\{x': \widehat{D}_{\alpha,m}(x,x')\leq h\}$, and $\mathcal{N}$ denotes the covering number.
Note that since $\widehat{\Gamma}(x,x')=\emptyset \Rightarrow \widehat{D}_{\alpha,m} = \infty$, we will always have $(x,x')\in\Psi$ if $x'\in S\cap S_{x,h}^{\widehat{D}_{\alpha,m}}$ (and, of course, the same applies when $x'\in S_m^*\cap S_{x,h/2}^{\widehat{D}_{\alpha,m}}$).
So we can apply Proposition \ref{thm:dbdest2} to give
\begin{align*}
\mathcal{G}_m \sup\limits_{x\in S_m^*} \sup\limits_{x'\in S\cap S_{x,h}^{\widehat{D}_{\alpha,m}}} D_\alpha(x,x')^{2\beta} \leq
\left[h \left(\frac{\lambda_0+\epsilon_m}{\lambda_0}\right)^\alpha \right]^{2\beta}
\end{align*}
and
\begin{align*}
\mathcal{G}_m \mathcal{N}\left(S_m^*,\widehat{D}_{\alpha,m},\frac{h}{2}\right) &\leq 
\mathcal{G}_m \mathcal{N}\left(S_m^*,d_{S_m^*},\frac{h(\lambda_0-\epsilon_m)^\alpha}{2}\right)
\end{align*}
where the $d_{S_m^*}$ distance is the length of the shortest path between two points restricted to $S_m^*$, as defined in the appendix.
Clearly $S_m^*$ has condition number at least $\tau_0-3\delta_m>0$.
If $S_m^*$ has exactly one connected component, then Proposition \ref{thm:geodesicdiam} combined with the assumption that $h > (2c_4/(\tau_0^{d-1}(\lambda_0-\epsilon_m)^\alpha)$ implies that any point in $S_m^*$ is a $h(\lambda_0-\epsilon_m)^\alpha/2$ covering, so
$$
\mathcal{N}\left(S_m^*,d_{S_m^*},\frac{h(\lambda_0-\epsilon_m)^\alpha}{2}\right)=1.
$$
Since $S_m^*$ can have at most $K$ connected components, we can repeat the same argument for each component and conclude that
$$
\mathcal{N}\left(S_m^*,d_{S_m^*},\frac{h(\lambda_0-\epsilon_m)^\alpha}{2}\right)\leq K.
$$
So,
\begin{align*}
&\mathbb{E}_{n,m}\left\{\int (\widehat{f}_{h,\alpha}(x)-f(x))^2 dP(x)\right\} \\&\leq
(M^2+\sigma^2)\left(\frac{1}{m}+3 c_3 2^d\Lambda_0 \frac{\delta_m}{\tau_0}\right) \\
&+ \left[h \left(\frac{\lambda_0+\epsilon_m}{\lambda_0}\right)^\alpha \right]^{2\beta}\\
&+ \frac{K(M^2/e+2\sigma^2)}{n}.
\end{align*}

\end{proof}

\subsubsection*{Acknowledgments}
This research is supported in part by AFOSR under grants 
FA9550-10-1-0382 and FA95500910373 and NSF under grants IIS-1116458 and DMS-0806009.


\bibliography{thesis}
\newpage
\section*{Appendix}

\noindent{\bf Results used in proof of Theorem~\ref{thm:SSL_UB}}\hspace{11pt}
In order to prove Theorem~\ref{thm:SSL_UB}, we characterize how the 
plug-in density-sensitive distance estimate $\hat D_\alpha$ behaves. For this, 
we start with a result about the density estimator.
\begin{thm}\label{thm:densityest}
If $m\geq m_0$, where $m_0\equiv m_0(\lambda_0,\Lambda_0)$ is a constant, then for all marginal densities $p$ of distributions in $\mathcal{P}_{XY}(\alpha)$, 
we have with probability $>1-1/m$,
\begin{align*}
\sup\limits_{x\in S\backslash \mathcal{R}_{\partial S}} |p(x)-\widehat{p}_m(x)|
\leq \epsilon_m \mbox{ and } \; \partial\widehat{S}\subset\mathcal{R}_{\partial S} 
\end{align*}
where $\epsilon_m=c_1 (\log m)^{-1/2}$ for constant $c_1\equiv c_1(K,C_2,d,\eta,\Lambda_0)$, $\widehat{S}=\{x:\widehat{p}_m(x)>0\}$, and
\begin{align*}
\mathcal{R}_{\partial S} = \left\{x:\inf\limits_{z\in\partial S} \|x-z\|_2 < \delta_m\right\}
\end{align*}
where $\delta_m=2c_2\sqrt{d}\left(\frac{\log^2m}{m}\right)^{\frac{1}{d}}$ for some constant $c_2>0$.
\end{thm}
\begin{proof}
Follows from Theorem 1 in \citet{ASingh:unlabeled} by noting that since the density estimate will be $0$ a.s. outside the boundary region, and we have $p\geq\lambda_0$ on $S$, for sufficiently large $m$ (i.e. small $\epsilon_m$), we must have $S\backslash \mathcal{R}_{\partial S}\subseteq\widehat{S}\subseteq S\cup \mathcal{R}_{\partial S}$.
\end{proof}


The following two propositions now characterize how the 
plug-in density-sensitive distance estimate $\hat D_\alpha$ behaves.
\begin{prop}\label{thm:dbdest1}
Assume $\sup\limits_{x\in S\backslash \mathcal{R}_{\partial S}}|\widehat{p}_m(x)-p(x)|\leq\epsilon_m$ and $\partial\widehat{S}\subset\mathcal{R}_{\partial S}$.
Let
\begin{align*}
\widetilde{D}_{\alpha,m}(x_1,x_2)=\inf\limits_{\gamma\in\widehat{\Gamma}(x_1,x_2)}\int\limits_{0}^{L(\gamma)} \frac{1}{p(\gamma(t))^\alpha}dt
\end{align*}
and $\Psi=\{(x_1,x_2): x_1,x_2\in \widehat{S}\backslash \widehat{\mathcal{R}}_{\partial S},\; \widehat{\Gamma}(x_1,x_2)\neq\emptyset\}$.
Then for any $(x_1,x_2)\in\Psi$,
\begin{align*}
\left(\frac{\lambda_0}{\lambda_0+\epsilon_m}\right)^\alpha \widetilde{D}_{\alpha,m}(x_1,x_2) & \leq \widehat{D}_{\alpha,m}(x_1,x_2) \\
& \hspace{-1.2cm} \leq \left(\frac{\lambda_0}{(\lambda_0-\epsilon_m)_+}\right)^\alpha  \widetilde{D}_{\alpha,m}(x_1,x_2).
\end{align*}
\end{prop}


\begin{proof}
Note that by the triangle inequality, $\mathcal{R}_{\partial S}\subseteq\widehat{\mathcal{R}}_{\partial S}$, so $\widehat{S}\backslash \widehat{\mathcal{R}}_{\partial S}\subseteq S\backslash \mathcal{R}_{\partial S}$
since $\tau_0 > 2\delta_m$ for $m$ large enough.
We see that if $(x_1,x_2)\in\Psi$, then $x$ and $y$ must be in the same connected component of $\widehat{S}\backslash \widehat{\mathcal{R}}_{\partial S}$, and, furthermore, all points along any path in $\widehat{\Gamma}(x_1,x_2)$ must also be in the same connected component.
For $(x_1,x_2)\in\Psi$,
\begin{align*}
&\widehat{D}_{\alpha,m}(x_1,x_2)\\
&=\inf\limits_{\gamma\in\widehat{\Gamma}(x_1,x_2)}\int\limits_{0}^{L(\gamma)} \frac{1}{p(\gamma(t))^\alpha} \frac{p(\gamma(t))^\alpha}{\widehat{p}_m(\gamma(t))^\alpha}dt\\
&\leq \inf\limits_{\gamma\in\widehat{\Gamma}(x_1,x_2)}  \left[\int\limits_{0}^{L(\gamma)} \frac{1}{p(\gamma(t))^\alpha} dt\right]\left[ \sup\limits_{t\in[0,L(\gamma)]} \left(\frac{p(\gamma(t))}{\widehat{p}_m(\gamma(t))}\right)^\alpha\right]\\
&\leq \sup\limits_{z\in S\backslash \mathcal{R}_{\partial S}} \left(\frac{p(z)}{\widehat{p}_m(z)}\right)^\alpha \widetilde{D}_{\alpha,m}(x_1,x_2)
\end{align*}
and
\begin{align*}
\sup\limits_{z\in S\backslash \mathcal{R}_{\partial S}} \left(\frac{p(z)}{\widehat{p}_m(z)}\right)^\alpha
&\leq \sup\limits_{z\in S\backslash \mathcal{R}_{\partial S}} \left(\frac{p(z)}{(p(z)-\epsilon_m)_+}\right)^\alpha \\
&\leq \left(\frac{\lambda_0}{(\lambda_0-\epsilon_m)_+}\right)^\alpha .
\end{align*}
So
\begin{align*}
\widehat{D}_{\alpha,m}(x_1,x_2)
&\leq \left(\frac{\lambda_0}{(\lambda_0-\epsilon_m)_+}\right)^\alpha \widetilde{D}_{\alpha,m}(x_1,x_2) .
\end{align*}

Similarly,

\begin{align*}
\widehat{D}_{\alpha,m}(x_1,x_2)
&\geq \inf\limits_{z\in S\backslash \mathcal{R}_{\partial S}} \left(\frac{p(z)}{p(z)+\epsilon_m}\right)^\alpha \widetilde{D}_{\alpha,m}(x_1,x_2) \\
&\geq \left(\frac{\lambda_0}{\lambda_0+\epsilon_m}\right)^\alpha \widetilde{D}_{\alpha,m}(x_1,x_2) .
\end{align*}

\end{proof}


Given a set $A\subseteq\mathbb{R}^d$, define
\begin{align*}
d_A(x_1,x_2)=\inf\limits_{\gamma\in\Gamma_A(x_1,x_2)} L(\gamma)
\end{align*}
where $\Gamma_A(x_1,x_2)=\{\gamma\in\Gamma(x_1,x_2): \forall t\in[0,L(\gamma)]\; \gamma(t)\in A\}$.

\begin{prop}\label{thm:dbdest2}
With the notation of Proposition \ref{thm:dbdest1}, for all $x_1,x_2$,
\begin{align*}
D_{\alpha,m}(x_1,x_2)\leq\widetilde{D}_{\alpha,m}(x_1,x_2).
\end{align*}
Assume $\sup\limits_{x\in S\backslash \mathcal{R}_{\partial S}}|\widehat{p}_m(x)-p(x)|\leq\epsilon_m$ and $\partial\widehat{S}\subset\mathcal{R}_{\partial S}$.
Then for any $(x_1,x_2)\in\Psi$,
\begin{align*}
\widetilde{D}_{\alpha,m}(x_1,x_2)\leq \frac{d_{S\backslash \widehat{\mathcal{R}}_{\partial S}}(x_1,x_2)}{\lambda^\alpha}
\end{align*}
and
\begin{align*}
\left(\frac{\lambda}{\lambda+\epsilon_m}\right)^\alpha D_\alpha(x_1,x_2) \leq \widehat{D}_{\alpha,m}(x_1,x_2) \leq \frac{d_{S_m^*}(x_1,x_2)}{(\lambda_0-\epsilon_m)_+^\alpha}
\end{align*}
where $S_m^*=\left\{x\in S:\inf\limits_{z\in\partial S} \|x-z\|_2 \geq 3\delta_m\right\}$.
\end{prop}

\begin{proof}
Since for any $x_1$ and $x_2$, $\widehat{\Gamma}(x_1,x_2)\subseteq\Gamma(x_1,x_2)$, clearly $D_{\alpha,m}(x_1,x_2)\leq\widetilde{D}_{\alpha,m}(x_1,x_2)$.
If $\partial\widehat{S}\subset\mathcal{R}_{\partial S}$, write
\begin{align*}
\widetilde{D}_{\alpha,m}(x_1,x_2)&=\inf\limits_{\gamma\in\widehat{\Gamma}(x_1,x_2)}\int\limits_{0}^{L(\gamma)} \frac{1}{p(\gamma(t))^\alpha}dt\\
&\leq \left[\sup\limits_{z\in \widehat{S}\backslash \widehat{\mathcal{R}}_{\partial S}}\frac{1}{p(z)^\alpha}\right] \left[\inf\limits_{\gamma\in\widehat{\Gamma}(x_1,x_2)}\int\limits_{0}^{L(\gamma)} dt\right]\\
&\leq \frac{1}{\lambda_0^\alpha} d_{\widehat{S}\backslash \widehat{\mathcal{R}}_{\partial S}}(x_1,x_2)\\
&\leq \frac{1}{\lambda_0^\alpha} d_{S_m^*}(x_1,x_2)
\end{align*}
since, by the triangle inequality, $S_m^*\subseteq\widehat{S}\backslash \widehat{\mathcal{R}}_{\partial S}$.
Applying Proposition \ref{thm:dbdest1}, the result follows.
\end{proof}


To prove Theorem~\ref{thm:SSL_UB}, we also 
need the following two results.

\begin{prop}\label{thm:condnumarea}
Let $\mathcal{X}$ be a compact subset of $\mathbb{R}^d$, and $T>0$.
Then for any $\tau\in(0,T)$, for all sets $S\subseteq\mathcal{X}$ with condition number at least $\tau$, $\Vol(\partial S)\leq c_3/\tau$ for some $c_3$ independent of $\tau$, where $\Vol$ is the $d-1$-dimensional volume.
\end{prop}
\begin{proof}
Let $\{z_i\}_{i=1}^N$ be a minimal Euclidean $\tau/2$-covering of $\partial S$, and $B_i=\{x:\|x-z_i\|_2\leq \tau/2\}$.
Let $T_i$ be the tangent plane to $\partial S$ at $z_i$.
Then using the argument made in the proof of Lemma 4 in \citet{genovese2010minimax},
\begin{align*}
\Vol(B_i\cap\partial S)&\leq C_1 \Vol(B_i\cap T_i) \frac{1}{\sqrt{1-(\tau/2)^2/\tau^2}}\\
&\leq C_2 \tau^{d-1}
\end{align*}
for some constants $C_1$ and $C_2$ independent of $\tau$.
Since $\mathcal{X}$ is compact,
\begin{align*}
\mathcal{N}(\partial S,\|\cdot\|_2,\tau/2)&\leq C \left(\frac{1}{\tau}\right)^d
\end{align*}
for some constant $C$ depending only on $\mathcal{X}$ and $T$, where $\mathcal{N}$ denotes the covering number (note that even though $\partial S$ is a $d-1$ dimensional set, we can't claim $\mathcal{N}(\partial S, \|\cdot\|_2, \tau)=O(\tau^{-(d-1)})$, since $\partial S$ can become space-filling as $\tau\rightarrow0$).
So
\begin{align*}
\Vol(\partial S) &\leq \sum\limits_{i=1}^N \Vol(B_i\cap\partial S)\\
&\leq C_2 \tau^{d-1} \mathcal{N}(\partial S,\|\cdot\|_2,\tau/2)\\
&\leq C_2C \tau^{-1}
\end{align*}
and the result follows with $c_3=C_2C$.

\end{proof}

\begin{prop}\label{thm:geodesicdiam}
Let $\mathcal{X}$ be a compact subset of $\mathbb{R}^d$, and $T>0$.
Then for any $\tau\in(0,T)$, for all compact, connected sets $S\subseteq\mathcal{X}$ with condition number at least $\tau$, $\sup\limits_{u,v\in S} d_{S} (u,v)\leq c_4 \tau^{1-d}$ for some $c_4$ independent of $\tau$.
\end{prop}
\begin{proof}
First consider the quantity
\begin{align*}
\sup\limits_{u,v\in\partial S} d_{S} (u,v).
\end{align*}
Since $\partial S\subseteq S$, clearly
\begin{align*}
\sup\limits_{u,v\in\partial S} d_{S} (u,v) \leq \sup\limits_{u,v\in\partial S} d_{\partial S} (u,v).
\end{align*}
Since $\partial S$ is closed, there must exist $u^*,v^*\in\partial S$ such that 
\begin{align*}
\sup\limits_{u,v\in\partial S} d_{\partial S} (u,v) = d_{\partial S}(u^*,v^*). 
\end{align*}
Let $\{z_i\}_{i=1}^{N}$ be a minimal $\tau$-covering of $\partial S$ in the $d_{\partial S}$ metric.
Let $\{\widetilde{z}_i\}_{i=1}^{\widetilde{N}}\subseteq\{z_i\}_{i=1}^{N}$ such that $d_{\partial S}(u^*,\widetilde{z}_1)\leq\tau$, $d_{\partial S}(v^*,\widetilde{z}_{\widetilde{N}})\leq\tau$, and for any $1\leq i\leq\widetilde{N}-1$, $d_{\partial S}(\widetilde{z}_{i},\widetilde{z}_{i+1})\leq2\tau$.
Then
\begin{align*}
d_{\partial S}(u^*,v^*) &\leq d_{\partial S}(u^*,\widetilde{z}_{1}) + d_{\partial S}(v^*,\widetilde{z}_{\widetilde{N}}) \\ &+ \sum\limits_{i=1}^{\widetilde{N}-1} d_{\partial S}(\widetilde{z}_{i},\widetilde{z}_{i+1})\\
&\leq 2\tau\widetilde{N}.
\end{align*}
So,
\begin{align*}
d_{\partial S}(u^*,v^*) &\leq 2 \tau \mathcal{N}(\partial S, d_{\partial S}, \tau).
\end{align*}
By Proposition 6.3 in \citet{niyogi2006finding} (or see Lemma 3 in \citet{genovese2010minimax}), if $x,y\in\partial S$ such that $\|x-y\|_2=a\leq\tau/2$, then $d_{\partial S}(x,y)\leq\tau-\tau\sqrt{1-(2a)/\tau}$.
In particular, if $\|x-y\|_2\leq\tau/2$, then $d_{\partial S}(x,y)\leq\tau$.
So any Euclidean $\tau/2$-covering of $\partial S$ is also a $\tau$-covering in the $d_{\partial S}$ metric.
Then we have
\begin{align*}
\sup\limits_{u,v\in\partial S} d_{S} (u,v) &\leq
d_{\partial S}(u^*,v^*) \\
&\leq 2 \tau \mathcal{N}(\partial S, d_{\partial S}, \tau)\\
&\leq 2 \tau \mathcal{N}(\partial S, \|\cdot\|_2, \tau/2)\\
&\leq C\tau \left(\frac{1}{\tau}\right)^d\\
&= C \tau^{1-d}
\end{align*}
for some constant $C$ depending only on $\mathcal{X}$ and $T$ (note that, as in the proof of \ref{thm:condnumarea}, even though $\partial S$ is a $d-1$ dimensional set, we can't claim $\mathcal{N}(\partial S, \|\cdot\|_2, \tau)=O(\tau^{-(d-1)})$, since $\partial S$ can become space-filling as $\tau\rightarrow0$).

Now let $u^\dag,v^\dag\in S$ such that
\begin{align*}
\sup\limits_{u,v\in S} d_S (u,v)=d_S(u^\dag,v^\dag)
\end{align*}
which must exist since $S$ is compact.
Let $u^\ddag,v^\ddag\in\partial S$ be the (not necessarily unique) projections of $u^\dag$ and $v^\dag$ onto $\partial S$.
Clearly the line segment connecting $u^\dag$ and $u^\ddag$ is fully contained in $S$, and the same applies to $v^\dag$ and $v^\ddag$.
So
\begin{align*}
d_S(u^\dag,v^\dag) &\leq d_S(u^\dag,u^\ddag) + d_S(u^\ddag,v^\ddag) + d_S(v^\ddag,v^\dag)\\
&\leq \|u^\dag-u^\ddag\|_2 + \|v^\dag-v^\ddag\|_2 + d_S(u^*,v^*)\\
&\leq 2\diam(\mathcal{X}) + C \tau^{1-d}
\end{align*}
and setting $c_4=2T^{d-1}\diam(\mathcal{X})+C$, the result follows.
\end{proof}

\noindent{\bf Proof of Theorem~\ref{thm:SL_LB}}\hspace{11pt}
The proof of Theorem~\ref{thm:SL_LB} is based on the following
result based on Assouad's Lemma (see e.g. \citet{Tsy2009}).

\begin{thm}
Let $\Omega = \{0,1\}^q$, the collection of binary vectors of length
$q\geq 1$. Let $\mathcal{P}_\Omega = \{P^\bomega,
\bomega\in\Omega\}$ be the corresponding collection of $2^q$
probability measures associated with each vector.  Also let
$\|P^{\bomega'} \wedge P^{\bomega}\|$ denote the affinity between
two distributions (i.e. $\|P^{\bomega'} \wedge P^{\bomega}\|=1-\sup\limits_{A}|P^{\bomega'}(A)-P^{\bomega}(A)|$,
where the supremum is over all measurable sets),
and $\rho(\cdot,\cdot)$ denotes the Hamming
distance between two binary vectors. For any semi-distance $d$
\begin{align*}
\inf_{\hat\bomega}\max_{\bomega\in\Omega} 
\mathbb{E}_{\bomega}[d^2(f^{\bomega},f^{\hat \bomega})] &\geq \frac{q}{8}
\left(\min_{\omega,\omega': \rho(\omega,\omega') \neq 0}
\frac{d^2(f^\omega,f^{\omega'})}{\rho(\omega,\omega')} \right)\\
&\times\left(\min_{\omega,\omega': \rho(\omega,\omega') =1 } \|P^{\bomega} \wedge P^{\bomega'}\|\right)
\end{align*}
\label{thm:minimax_bound}
\end{thm}

We now prove Theorem~\ref{thm:SL_LB}.

\begin{proof}

\textbf{Construction:}

Let $l=\lfloor c_0 n^{1/(d-1)}\rfloor$ with $c_0>1$ a constant, $q=l^{d-1}$, $\Omega=\{0,1\}^q$ and $\epsilon=\frac{1}{l+2}$. 
For $i\in\{1,...,l\}$, let $a_i=\frac{i+0.5}{l+2}$.
For $\vec{i}\in\{1,...,l\}^{d-1}$, let $v_{\vec{i}}=(a_{\vec{i}_1},...,a_{\vec{i}_{d-1}})$. 
Define $g:\mathbb{R}^{d-1}\rightarrow\mathbb{R}$ as $g(\tilde{x}) = $
\begin{align*}
\left\{
\begin{array}{rcl}
r+\sqrt{\left(\frac{1}{2}-r\right)^2-\|\tilde{x}\|_2^2} & \mbox{for} & \|\tilde{x}\|_2<\frac{1}{2}-r \\
r-\sqrt{r^2-\left(\frac{1}{2}-\|\tilde{x}\|_2\right)^2} & \mbox{for} & \frac{1}{2}-r\leq\|\tilde{x}\|_2<\frac{1}{2} \\
0 & \mbox{o.w.} & 
\end{array}\right.
\end{align*}
for $\tilde{x}\in\mathbb{R}^{d-1}$, where $r\in(0,1/4)$, to be specified later.
Let $B=\{(\tilde{x},x_d)\in[-0.5,0.5]^{d-1}\times[0,1]:x_d\leq g(\tilde{x})\}$.
For $\vec{i}\in\{1,...,l\}^{d-1}$, let $\underline{B}_{\vec{i}}=\{(\tilde{x},x_d)\in\mathbb{R}^{d-1}\times\mathbb{R}:((\tilde{x}-v_{\vec{i}})/\epsilon,x_d-1/8)\in B \}$ and $\overline{B}_{\vec{i}}=\{(\tilde{x},x_d)\in\mathbb{R}^{d-1}\times\mathbb{R}:((\tilde{x}-v_{\vec{i}})/\epsilon,x_d-(1/8+r))\in B \}$.
Let $\underline{S}=\{x\in\mathbb{R}^d:\exists x'=(\tilde{x}',x_d')\in[\epsilon,1-\epsilon]^{d-1}\times[\epsilon,\frac{1}{8}-\epsilon] \text{ s.t. } \|x-x'\|_2\leq\epsilon\}$ and $\overline{S}=\{x\in\mathbb{R}^d:\exists x'=(\tilde{x}',x_d')\in[\epsilon,1-\epsilon]^{d-1}\times[\frac{1}{8}+r+\epsilon,1-\epsilon] \text{ s.t. } \|x-x'\|_2\leq\epsilon\}$.
For any $\Gamma\subseteq\{1,...,l\}^{d-1}$, let $\underline{S}_\Gamma=\underline{S}\cup\left(\bigcup\limits_{\vec{i}\in\Gamma}\underline{B}_{\vec{i}}\right)$ and $\overline{S}_\Gamma=\overline{S}\backslash\left(\bigcup\limits_{\vec{i}\in\Gamma}\overline{B}_{\vec{i}}\right)$.
Let $\vec{\Gamma}$ be an arbitrary ordering of $\{1,...,l\}^{d-1}$.
Given $\bomega\in\Omega$, let $\Gamma(\bomega)=\{\vec{\Gamma}_i:\bomega_i=1\}$, and let $\underline{S}^\bomega=\underline{S}_{\Gamma(\bomega)}$, $\overline{S}^\bomega=\overline{S}_{\Gamma(\bomega)}$, and $S^\bomega=\underline{S}^\bomega\cup\overline{S}^\bomega$.

Let $p^\bomega(x)=\frac{I_{S^\bomega}(x)}{\leb(S^\bomega)}$, $f^\bomega(x)=MI_{\underline{S}^\bomega}(x)-MI_{\overline{S}^\bomega}(x)$, and $p(y|x)^\bomega=\delta(y-f^\bomega(x))$, where $\delta(\cdot)$ is the Dirac delta (we could also use a conditional distribution that is absolutely continuous with respect to Lebesgue measure; the result would be the same).
Finally, let $P^\bomega$ denote the measure on $\mathbb{R}^{d+1}$ defined by $p^\bomega(x)$ and $p^\bomega(y|x)$, and $P^\bomega_n$ the corresponding product measure.

\textbf{Proof of $\Omega(1)$ rate:}

Note that $\leb(\underline{B}_{\vec{i}})=\leb(\overline{B}_{\vec{i}})$, and so for any $\bomega,\bomega'$, $\leb(S^\bomega)=\leb(S^{\bomega'})=\leb(\underline{S})+\leb(\overline{S})$.
Let $\lambda=1/(\leb(\underline{S})+\leb(\overline{S}))$, i.e. $\lambda=1/\leb(S^\bomega)$ for any $\bomega$.

Let $\bomega,\bomega'\in\Omega$ such that $\rho(\bomega,\bomega')=1$ (where $\rho$ denotes the hamming distance), and WLOG assume $\bomega_i=0$ and $\bomega'_i=1$.
Also denote $\vec{i}=\vec{\Gamma}_i$.
Then the L1 distance between $P^\bomega$ and $P^{\bomega'}$ is
\begin{align*}
&d_1(P^\bomega,P^{\bomega'})\\&=\int\limits_{\mathbb{R}^d}\int\limits_{\mathbb{R}} |p^\bomega(x)p^\bomega(y|x)-p^{\bomega'}(x)p^{\bomega'}(y|x)|dydx\\
&=\int\limits_{\underline{S}^\bomega\cup\overline{S}^{\bomega'}}\int\limits_{\mathbb{R}} |\lambda p^\bomega(y|x)-\lambda p^{\bomega'}(y|x)|dydx\\
&+\int\limits_{\overline{B}_{\vec{i}}\backslash\underline{B}_{\vec{i}}}\int\limits_{\mathbb{R}} \lambda p^\bomega(y|x)dydx
+\int\limits_{\underline{B}_{\vec{i}}\backslash\overline{B}_{\vec{i}}}\int\limits_{\mathbb{R}} \lambda p^{\bomega'}(y|x)dydx\\
&+\int\limits_{\underline{B}_{\vec{i}}\cap\overline{B}_{\vec{i}}}\int\limits_{\mathbb{R}} |\lambda p^\bomega(y|x)-\lambda p^{\bomega'}(y|x)|dydx\\
&= 0 + \lambda \leb(\overline{B}_{\vec{i}}\backslash\underline{B}_{\vec{i}}) + \lambda \leb(\underline{B}_{\vec{i}}\backslash\overline{B}_{\vec{i}}) \\ &\;\;\;+ 2 \lambda \leb(\underline{B}_{\vec{i}}\cap\overline{B}_{\vec{i}})\\
&=\lambda (\leb(\underline{B}_{\vec{i}})+\leb(\overline{B}_{\vec{i}}))\\
&=2\lambda\epsilon^{d-1}\leb(B)
\end{align*}
where in the first step we have used the fact that $x\notin S^\bomega\cup S^{\bomega'}\Rightarrow p^\bomega(x)=p^{\bomega'}(x)=0$, and divided $S^\bomega\cup S^{\bomega'}$ into four non-intersecting components.
Then we can bound the affinity of the product measures $P^\bomega_n$ and $P^{\bomega'}_n$ for $\rho(\bomega,\bomega')=1$ as
\begin{align*}
\|P^\bomega_n\wedge P^{\bomega'}_n\|&\geq(1-d_1(P^\bomega,P^{\bomega'})/2)^n\\
&=(1-\lambda\epsilon^{d-1}\leb(B))^n.
\end{align*}

For any $\bomega\neq\bomega'$, denoting as $\bomega\wedge\bomega'$ the logical and of $\bomega$ and $\bomega'$, we have, for arbitrary $\vec{j}\in\{1,...,l\}^{d-1}$,
\begin{align*}
&d^2(f^\bomega,f^{\bomega'}) \\
&=\sum\limits_{\vec{i}\in\Gamma(\bomega\wedge\bomega')}\int\limits_{\underline{B}_{\vec{i}}\Delta\overline{B}_{\vec{i}}} M^2 dx+\int\limits_{\underline{B}_{\vec{i}}\cap\overline{B}_{\vec{i}}} 4M^2 dx\\
&=\rho(\bomega,\bomega')(M^2\leb(\underline{B}_{\vec{j}}\Delta\overline{B}_{\vec{j}})+4M^2\leb(\underline{B}_{\vec{j}}\cap\overline{B}_{\vec{j}}))\\
&=2\rho(\bomega,\bomega')M^2(\leb(\underline{B}_{\vec{j}})+\leb(\underline{B}_{\vec{j}}\cap\overline{B}_{\vec{j}}))\\
&=2\rho(\bomega,\bomega')M^2\epsilon^{d-1}(\leb(B)+\leb(B_r))
\end{align*}
where we define $B_r=\{x\in B:x-(0,...,0,r)\in B\}$.
Then by Theorem \ref{thm:minimax_bound},
\begin{align*}
&\inf_{\hat\bomega}\max_{\bomega\in\Omega} 
\mathbb{E}_{\bomega}[d^2(f^{\bomega},f^{\hat \bomega})] 
\\&\geq \frac{M^2 (l\epsilon)^{d-1}}{4}(\leb(B)+\leb(B_r))(1-\lambda\epsilon^{d-1} \leb(B))^n .
\end{align*}
Also we have
\begin{align*}
&\frac{1}{\lambda}\int(f^\bomega(x)-f^{\bomega'}(x))^2p^\bomega(x)dx = \int\limits_{S^\bomega} (f^\bomega(x)-f^{\bomega'}(x))^2 dx\\
&= \int\limits (f^\bomega(x)-f^{\bomega'}(x))^2 dx - \int\limits_{S^{\bomega'}\backslash S^\bomega} (f^\bomega(x)-f^{\bomega'}(x))^2 dx\\
&= d^2(f^\bomega,f^{\bomega'}) - M^2 \vol(S^{\bomega'}\backslash S^\bomega)\\
&\geq d^2(f^\bomega,f^{\bomega'}) - M^2q \epsilon^{d-1} \leb(B\backslash B_r)\\
&= d^2(f^\bomega,f^{\bomega'}) - M^2 \left(\frac{l}{l+2}\right)^{d-1}(\leb(B)-\leb(B_r)).
\end{align*}

Since $\lambda>1$,
\begin{align*}
&\inf_{\widehat{f}}\sup_{(p,f)\in\mathcal{P}_{XY}(\alpha)} \mathbb{E}_n\int(\widehat{f}(x)-f(x))^2 dP(x)\\
&\geq \inf_{\hat\bomega}\max_{\bomega\in\Omega} \mathbb{E}_{\bomega}\int(f^{\hat\bomega}(x)-f^\bomega(x))^2p^\bomega(x)dx \\
&\geq  \frac{M^2 (l\epsilon)^{d-1}}{4}(\leb(B)+\leb(B_r))(1-\lambda\epsilon^{d-1} \leb(B))^n \\
&-M^2 \left(\frac{l}{l+2}\right)^{d-1}(\leb(B)-\leb(B_r)).
\end{align*}

Assume $n\geq 2^d$.
Then $l\geq 2$ and
\begin{align*}
\left(\frac{l}{l+2}\right)^{d-1}\geq \frac{1}{2^{d-1}}.
\end{align*}
Clearly $\leb(B)\leq\frac{1}{2}$.
Let $c_0\geq3$.
Then $\epsilon\leq1/8$ and $\lambda\leq(1-2\epsilon)^{-(d-1)}(1-4\epsilon - r)^{-1}\leq 2^{d+1}$, so
\begin{align*}
(1-\lambda\epsilon^{d-1} \leb(B))^n\geq \left(1-\frac{2^d}{c_0^{d-1}n} \right)^n\rightarrow e^{-2^d/c_0^{d-1}}.
\end{align*}
So if we let $c_0>(2^d/\log(5/4))^{1/(d-1)}$, then $e^{-2^d/c_0^{d-1}}>4/5$ and for sufficiently large $n$ we will have $(1-\lambda\epsilon^{d-1} \leb(B))^n\geq4/5$.
Hence,
\begin{align*}
\inf_{\widehat{f}}\sup_{(p,f)\in\mathcal{P}_{XY}(\alpha)} \mathbb{E}_n\int(\widehat{f}(x)-f(x))^2 dP(x)\\
\geq  \frac{M^2 }{5\cdot2^{d-2}}\left(\leb(B_r) - 2\leb(B\backslash B_r)\right).
\end{align*}

Since
\begin{align*}
\leb(B_r) &= \frac{1}{2}\left(\frac{1}{2}-r\right)^d\frac{\pi^{d/2}}{\Gamma(d/2+1)}
\end{align*}
and
\begin{align*}
\leb(B\backslash B_r) \leq r\frac{\pi^{(d-1)/2}}{2^{d-1}\Gamma((d-1)/2+1)}
\end{align*}
where $\Gamma$ is the gamma function, then
\begin{align*}
&\leb(B_r) - 2\leb(B\backslash B_r) \\&\geq
\frac{1}{2}\left(\frac{1}{2}-r\right)^d\frac{\pi^{d/2}}{\Gamma(d/2+1)} - r\frac{\pi^{(d-1)/2}}{2^{d-2}\Gamma((d-1)/2+1)} \\
&\geq
\frac{\pi^{d/2}}{2^d\Gamma\left(\frac{d+1}{2}\right)d}\left[(1-2r)^d - \frac{4d}{\sqrt{\pi}}r\right].
\end{align*}
Now let $r$ be such that
\begin{align*}
(1-2r)^d - \frac{4d}{\sqrt{\pi}}r = \frac{1}{2}
\end{align*}
(it is easy to see that this can be satisfied by some $r\in(0,1/4)$ for any $d\geq1$).
So we have
\begin{align*}
\inf_{\widehat{f}}\sup_{(p,f)\in\mathcal{P}_{XY}(\alpha)} \mathbb{E}_n\int(\widehat{f}(x)-f(x))^2 dP(x)\\
\geq  \frac{M^2 \pi^{d/2}}{5\cdot 2^{2d-1}\Gamma\left(\frac{d+1}{2}\right)d}.
\end{align*}

\textbf{Verifying condition number:}

Let $\tau(A)$ be the condition number of a set $A$.
Then for arbitrary $\bomega$,
\begin{align*}
\tau(S^\bomega) = \min\left\{\tau(\underline{S}^\bomega),\tau(\overline{S}^\bomega),\frac{1}{2}\inf\limits_{u\in\underline{S}^\bomega}\inf\limits_{v\in\overline{S}^\bomega}\|u-v\|_2\right\}.
\end{align*}
Due to the shape of the function $g$, for arbitrary $\vec{i}\in\{1,...,l\}^{d-1}$ we have
\begin{align*}
\tau(\underline{S}^\bomega) &\geq \min\left\{\tau(\partial\underline{S}),\tau(\partial\underline{B}_{\vec{i}}\backslash\partial\underline{S})\right\}
\end{align*}

By definition of $\underline{S}$ it is easy to see that $\tau(\partial\underline{S})=\epsilon$.
Also
\begin{align*}
&\tau(\partial\underline{B}_{\vec{i}}\backslash\partial\underline{S})\\ &= \tau(\{(\tilde{x},x_d)\in[-\epsilon/2,\epsilon/2]^{d-1}\times[0,1]:x_d=g(\tilde{x}/\epsilon)\})\\
&= \epsilon \tau(\{(\tilde{x},x_d)\in[-1/2,1/2]^{d-1}\times[0,1]:x_d=g(\tilde{x})\})\\
&= \epsilon r.
\end{align*}
Since $r<1$, we have $\tau(\underline{S}^\bomega) \geq r\epsilon$, and similarly $\tau(\overline{S}^\bomega) \geq r\epsilon$.
Now,
\begin{align*}
&\frac{1}{2}\inf\limits_{u\in\underline{S}^\bomega}\inf\limits_{v\in\overline{S}^\bomega}\|u-v\|_2 \\&\geq
\frac{\epsilon}{2} \inf\limits_{u,v\in[-0.5,0.5]^{d-1}} \|(u,g(u))-(v,g(v)+r)\|_2\\
&=\frac{\epsilon}{2} \left(\sqrt{\frac{1}{4}+r^2}-\frac{1}{2}\right)
\end{align*}
which is smaller than $\epsilon r$, so for $n$ sufficiently large,
\begin{align*}
\tau(S^\bomega)&\geq\frac{\epsilon}{2} \left(\sqrt{\frac{1}{4}+r^2}-\frac{1}{2}\right)\\
&\geq\frac{1}{2(c_0n^{1/(d-1)}+2)} \left(\sqrt{\frac{1}{4}+r^2}-\frac{1}{2}\right)\\
&\geq n^{-\frac{1}{d-1}}\frac{1}{2(c_0+1)} \left(\sqrt{\frac{1}{4}+r^2}-\frac{1}{2}\right)
\end{align*}
which completes the proof.

\end{proof}



\noindent{\bf Proof of Theorem~\ref{thm:crossval}}\hspace{11pt}
%
%
%
%
%
First, we derive a general concentration of $\hat \cE(f)$ around $\cE(f)$
where 
$\hat \cE(f) = \hat R(f) - \hat R(f^*) = -\frac1{n}\sum^n_{i=1}U_i,$
and $U_i = -(Y_i-f(X_i))^2+(Y_i-f^*(X_i))^2$. 

If the variables $U_i$ satisfy the following moment condition:
$$
\E[|U_i - \E[U_i]|^k] \leq \frac{var(U_i)}{2}k! r^{k-2}
$$
for some $r>0$, then the Craig-Bernstein (CB) inequality (Craig 1933)
states that with probability $>1-\delta$,
$$
\frac1{n}\sum^n_{i=1} (U_i - \E[U_i]) \leq \frac{\log(1/\delta)}{nt} + \frac{t \ var(U_i)}{2(1-c)}
$$
for $0\leq tr \leq c <1$. The moment conditions are satisfied by
bounded random variables as well as Gaussian random variables (see e.g.
\citet{Randproj_JHaupt}).

To apply this inequality, we first show that var($U_i) \leq 4(M^2+\sigma^2)\cE(f)$ since
$Y_i = f(X_i) + \epsilon_i$ with $\epsilon_i \stackrel{i.i.d}{\sim} {\cal N}(0,\sigma^2)$.
Also, we assume that $|f(x)|$, $|\hat f(x)| \leq M$, where $M>0$ is a constant.
\begin{eqnarray*}
&& \hspace{-0.8cm}var(U_i) \leq \E[U_i^2] \\
&& \hspace{-0.8cm}= \E[(-(Y_i-f(X_i))^2+(Y_i-f^*(X_i))^2)^2]\\
&& \hspace{-0.8cm}= \E[(-(f^*(X_i)+\epsilon_i-f(X_i))^2+(\epsilon_i)^2)^2]\\
&& \hspace{-0.8cm}= \E[(-(f^*(X_i)-f(X_i))^2-2\epsilon_i(f^*(X_i) - f(X_i)))^2]\\
&& \hspace{-0.8cm}\leq 4 M^2 \cE(f) + 4 \sigma^2 \cE(f) = 4(M^2+\sigma^2)\cE(f)
\end{eqnarray*}

Therefore using CB inequality we get, with probability 
$>1-\delta$,
$$
\cE(f) - \hat \cE(f) \leq \frac{\log(1/\delta)}{nt} + \frac{t \ 2(M^2+\sigma^2)\cE(f)}{(1-c)}
$$
Now set $c = tr = 8t(M^2+\sigma^2)/15$ and let $t <  15/(38(M^2+\sigma^2))$. With this choice,
$c < 1$ and define 
$$
a = \frac{t 2(M^2+\sigma^2)}{(1-c)} <1.
$$
Then, using $a$ and rearranging terms, with probability $>1-\delta$,
$$
(1-a)\cE(f) - \hat \cE(f) \leq \frac{\log(1/\delta)}{nt} 
$$
where $t <  15/(38(M^2+\sigma^2))$.

Then, using the previous 
concentration result, and taking union bound over all $f\in {\cal F}$, we have with probability $>1-\delta$,
$$
\cE(f) \leq \frac1{1-a}\left[\hat \cE^V(f) + \frac{\log(|{\cal F}|/\delta)}{nt} \right].
$$
Now consider
\begin{eqnarray*}
\cE(\hat f_{\hat \alpha,\hat h}) &&\hspace{-0.6cm}=  R(\hat f_{\hat \alpha,\hat h}) - R(f^*) \\
&&\hspace{-1cm}\leq 
\frac1{1-a}\left[\hat R^V(\hat f_{\hat \alpha,\hat h}) - \hat R^V(f^*) + \frac{\log(|{\cal F}|/\delta)}{nt} \right]\\
&&\hspace{-1cm}\leq  \frac1{1-a}\left[\hat R^V(f) - \hat R^V(f^*) + \frac{\log(|{\cal F}|/\delta)}{nt} \right]
\end{eqnarray*}
Taking expectation with respect to validation dataset,
\begin{eqnarray*}
\E_V[\cE(\hat f_{\hat \alpha,\hat h})] 
&\leq &\frac1{1-a}\left[R(f) - R(f^*) + \frac{\log(|{\cal F}|/\delta)}{nt} \right] \\
&&+ 4\delta M^2. 
\end{eqnarray*}
Now taking expectation with respect to training dataset,
\begin{eqnarray*}
\E_{TV}[\cE(\hat f_{\hat \alpha,\hat h})] 
&\leq& \frac1{1-a}\left[\E_T[R(f) - R(f^*)] \right.\\
&& +\left.\frac{\log(|{\cal F}|/\delta)}{nt} \right] + 
4 \delta M^2.
\end{eqnarray*}
Since this holds for all $f \in {\cal F}$, we get: 
\begin{eqnarray*}
\E_{TV}[\cE(\hat f_{\hat \alpha,\hat h})] 
& \leq & \frac1{1-a}\left[\min_{f\in{\cal F}}\E_T[\cE(f)] + 
\frac{\log(|{\cal F}|/\delta)}{nt} \right]\\
& &+  4 \delta M^2.
\end{eqnarray*}
The result follows since 
${\cal F} = \{\hat f^T_{\alpha,h}\}_{\alpha \in {\cal A}, h\in {\cal H}}$ and 
$|{\cal F}| = |{\cal A}| |{\cal H}|$.


\flushright{$\Box$}

%
%
%

\end{document}